\def \beq{\begin{equation}}
\def \eq{\end{equation}}
\def \berr{\begin{eqnarray}}
\def \err{\end{eqnarray}}
\def \a{\alpha}
\def \b{\beta}

\def \eps{\varepsilon}

\def \om{\omega}
\def \l{\lambda}
\def \dl{\partial}

\def \R{{\cal R}}

\def\Sq{S_q^{N-1}}

\def \({\left(}
\def \){\right)}
\def \<{\langle}
\def \>{\rangle}

\def \und#1{\underline{#1}}

\def\tens{\mathop{\otimes}}

\renewcommand \d{d^{N}\!}
\newcommand \reals{I \! \! R}

\def\lform{\hbox{$\sqcup$}\llap{\hbox{$\sqcap$}}}

\def\proof{\goodbreak\noindent{\bf Proof\quad}}
\def\endproof{{\ $\lform$}\bigskip }

\newtheorem{prop}{Proposition}[section]
\newtheorem{theorem}[prop]{Theorem}
\newtheorem{lemma}[prop]{Lemma}

\documentstyle[12pt]{article}

\advance \topmargin by -\headheight
\advance \topmargin by -\headsep
\evensidemargin \oddsidemargin
\begin{document}
\begin{titlepage}
\begin{center}
August 3, 1995
        \hfill  LBL-37431 \\
          \hfill    UCB-PTH-95/21 \\
\vskip .5in

{\large Integration on quantum Euclidean space and sphere in $N$ dimensions}
\footnote[1]{This work was supported in part by the Director, Office of
Energy Research, Office of High Energy and Nuclear Physics, Division of
High Energy Physics of the U.S. Department of Energy under Contract
DE-AC03-76SF00098 and in part by the National Science Foundation under
grant PHY90-21139.}

\vskip .3in

Harold Steinacker   \footnote[2]{email: hsteinac@physics.berkeley.edu}

\vskip .5in
{\em
Department of Physics,
University of California \\ and\\
   Theoretical Physics Group, \\ Lawrence Berkeley Laboratory\\
      1 Cyclotron Road, Berkeley, CA 94720}
\end{center}

\vskip .3in

\begin{abstract}
Invariant integrals of functions and forms over $q$ - deformed Euclidean space
and spheres in $N$ dimensions are defined and shown to be positive definite,
compatible with the star - structure and to satisfy a cyclic property involving
the $D$ - matrix of $SO_q(N)$. The definition is more general than the
Gaussian integral known so far. Stokes theorem is proved with and without
spherical boundary terms, as well as on the sphere.

\end{abstract}
\end{titlepage}
\renewcommand{\thepage}{\roman{page}}
\setcounter{page}{2}
\mbox{ }

\vskip 1in

\begin{center}
{\bf Disclaimer}
\end{center}

\vskip .2in

\begin{scriptsize}
\begin{quotation}
This document was prepared as an account of work sponsored by the United
States Government.  Neither the United States Government nor any agency
thereof, nor The Regents of the University of California, nor any of
their
employees, makes any warranty, express or implied, or assumes any legal
liability or responsibility for the accuracy, completeness, or
usefulness
of any information, apparatus, product, or process disclosed, or
represents
that its use would not infringe privately owned rights.  Reference
herein
to any specific commercial products process, or service by its trade
name,
trademark, manufacturer, or otherwise, does not necessarily
constitute or
imply its endorsement, recommendation, or favoring by the United States
Government or any agency thereof, or The Regents of the University of
California.  The views and opinions of authors expressed herein do not
necessarily state or reflect those of the United States
Government or any
agency thereof of The Regents of the University of California and shall
not be used for advertising or product endorsement purposes.
\end{quotation}
\end{scriptsize}

\vskip 2in

\begin{center}
\begin{small}
{\it Lawrence Berkeley Laboratory is an equal opportunity employer.}
\end{small}
\end{center}
\newpage
\renewcommand{\thepage}{\arabic{page}}
\setcounter{page}{1}

\section{Introduction}

In recent years, there has been much interest in formulating physics and in
particular field theory on quantized,
i.e. noncommutative spacetime. One of the motivations is that if there are
no more "points" in spacetime, such a theory should be well - behaved
in the UV. The concept of integration on such a space can certainly be
expected to be an essential ingredient.
In the simplest case of the quantum plane, such an integral was first
introduced by Wess and Zumino \cite{WZ}; see also \cite{CZ}.
In the presumably more physical case of quantum  Euclidean space \cite{FRT},
the Gaussian integration method was proposed by a number of authors
\cite{fiore,KM}. However,
it is very tedious to calculate except in the simplest cases and
its properties have never been investigated thoroughly; in fact,
its domain of definition turns out to be rather small.

In this paper, we will give a different definition based on spherical
integration in $N$ dimensions and investigate its properties in detail.
Although this idea has already appeared in the literature \cite{weich},
it has not been developed very far.
It turns out to be both simpler and more general than the Gaussian integral.
We first show that there is a unique invariant
integral over the quantum Euclidean sphere, and prove that it is positive
definite and satisfies a cyclic property involving the $D$ - matrix of
$SO_q(N)$. The integral over quantum Euclidean space is then defined by radial
integration, both for functions and $N$
forms. It turns out not to be
determined uniquely by rotation and translation invariance (=Stokes theorem)
alone; it is unique after e.g. imposing a general scaling law.
It is positive definite as well and thus allows to define a Hilbertspace of
square - integrable functions, and satisfies the same cyclic property.
The cyclic property also holds for the integral of $N$ and $N-1$ - forms over
spheres, which leads to a simple, truly noncommutative proof
of Stokes theorem on Euclidean space with and without spherical boundary
terms, as well as on the sphere.
These proofs only work for $q\neq 1$,
nevertheless they reduce to the classical Stokes theorem for $q\rightarrow 1$.
This shows the power of noncommutative geometry.
Obviously one would like to use this integral to define actions for
field theories on such noncommutative spaces; this is work in progress.

Although only the case of quantum Euclidean space is considered, the general
approach is clearly applicable to e.g. quantum Minkowski space as well.

\section{Integral on the quantum sphere $\Sq$}

To establish the notation, we briefly summarize the definitions used in this
paper, following Faddeev, Reshetikhin and Takhtadjan \cite{FRT}.

The (function algebra on the) quantum orthogonal group $O_q(N)$
(which is called $SO_q(N)$ in \cite{FRT})
is the algebra generated by $A^i_j$ modulo the relations
\berr
\hat R^{ik}_{mn} A^m_j A^n_l &=& A^i_n A^k_m \hat R^{nm}_{jl}, \label{RTT}\\
           g_{ij}A^i_k A^j_l &=& g_{kl}.
\err
$SO_q(N)$ is obtained by further imposing
\beq
A^{i_1}_{j_1} .... A^{i_N}_{j_N} \eps_q^{j_1 .... j_N} =
  \eps_q^{i_1 .... i_N}
\eq
using the fact that the quantum determinant is central, see e.g. \cite{meyer}.

The $\hat R$ - matrix decomposes into 3 projectors
$R^{ij}_{kl} = (q P^+ - q^{-1} P^- +q^{1-N} P^0)^{ij}_{kl}$.
The metric is determined by
$(P^0)^{ij}_{kl} =\frac {q^2-1}{(q^N-1)(q^{2-N}+1)} g^{ij} g_{kl}$,
where $g_{ik}g^{kj} = \delta_i^j$.
In this paper, we assume $q$ is real and positive. Then there is a star -
structure (involution)
\beq
\overline{A^i_j} = g^{jm} A^l_m g_{li}  \label{conj_A}
\eq
so that we really have $(S)O_q(N,\reals)$, and the antipode can be written as
\beq
 S(A^i_j) = \overline{A^j_i}.
\eq
Quantum Euclidean space \cite{FRT} is generated by $x^i$
with commutation relations
\beq
(P^-)^{ij}_{kl} x^k x^l =0,
\eq
and the center is generated by 1 and $r^2 = g_{ij} x^i x^j$.
The associated differentials satisfy
$(P^+)^{ij}_{kl} dx^k dx^l =0$ and  $g_{ij} dx^i dx^j = 0$, i.e.
\beq
dx^i dx^j = -q \hat R^{ij}_{kl} dx^k dx^l.
\eq
The epsilon - tensor is then determined by the unique top - (N-) form
\beq
dx^{i_1} ... dx^{i_N} = \eps_q^{i_1 ... i_N} dx^1 ... dx^N \equiv
\eps_q^{i_1 ... i_N} \d x.
\eq
The above relations are preserved under the coaction of $(S)O_q(N)$
\beq   \label{left_coaction}
\Delta(x^i) = A^i_j \tens x^j \equiv x^i_{(1)} \tens x^i_{(2)},
\eq
in Sweedler - notation. The involution $\overline{x^i} = x^j g_{ji}$ is
compatible with the left coaction of $(S)O_q(N,\reals)$.
One can also introduce derivatives which satisfy
\beq
(P^-)^{ij}_{kl} \dl^k \dl^l =0,
\eq
\beq
\dl^i x^j = g^{ij} + q (\hat R^{-1})^{ij}_{kl} x^k \dl^l,
\eq
and
\beq
\dl^i dx^j = q^{-1}\hat R^{ij}_{kl} dx^k \dl^l, \quad  x^i dx^j =
q \hat  R^{ij}_{kl} dx^k x^l.
\eq
This represents one possible choice. For more details, see e.g. \cite{OZ}.
Finally, the quantum sphere $\Sq$ is generated by $t^i = x^i/r$, which satisfy
$g_{ij} t^i t^j = 1$.

We first define a (complex - valued) integral $<f(t)>_t$ of a function $f(t)$
over $\Sq$.
It should certainly be invariant under $O_q(N)$,
which means
\beq
A^{i_1}_{j_1} ... A^{i_n}_{j_n} <t^{j_1} ... t^{j_n}>_t =
<t^{i_1} ... t^{i_n}>_t .
\eq
Of course, it has to satisfy
\beq
g_{i_l i_{l+1}} <t^{i_1} ... t^{i_n}>_t  =
<t^{i_1} ...t^{i_{l-1}}t^{i_{l+2}} ... t^{i_n}>_t \quad {\rm and} \quad
(P^-)^{i_l i_{l+1}}_{j_l j_{l+1}} <t^{j_1} ... t^{j_n}>_t = 0
\eq
We require one more property, namely that $I^{i_1 ... i_n} \equiv
<t^{i_1} ... t^{i_n}>_t$ is analytic in $(q-1)$, i.e. its Laurent series
in $(q-1)$ has no negative terms (we can then assume that the classical limit
$q=1$ is nonzero). These properties in fact determine the spherical integral
uniquely: for $n$ odd, one should define $<t^{i_1} ... t^{i_n}>_t =0$, and

\begin{prop}\label{spher_int}
For even $n$, there exists  (up to normalization) one and only one tensor
$I^{i_1 ... i_n} = I^{i_1 ... i_n}(q)$ analytic in $(q-1)$ which is invariant
under $O_q(N)$
\beq
A^{i_1}_{j_1} ... A^{i_n}_{j_n} I^{j_1 ... j_n} = I^{i_1 ... i_n}
\label{I_inv}
\eq
and symmetric,
\beq
(P^-)^{i_l i_{l+1}}_{j_l j_{l+1}} I^{j_1 ... j_{n}} = 0
\label{I_symm}
\eq
for any l. It can be normalized such that
\beq
g_{i_l i_{l+1}} I^{i_1 ... i_{n}} = I^{i_1 ... i_{l-1} i_{l+2} ... i_n}
\label{I_norm}
\eq
for any l. $I^{i j} \propto g^{i j}$.

An explicit form is e.g.
$I^{i_1 ... i_{n}} = \l_n (\Delta^{n/2} x^{i_1} ... x^{i_{n}})$, where
$\Delta=g_{ij} \dl^i \dl^j$ is the Laplacian
(in either of the 2 possible calculi),
and $\l_n$ is analytic in $(q-1)$.
For $q=1$, they reduce to  the  classical symmetric invariant tensors.
\end{prop}
\proof
The proof is by induction on $n$. For $n=2$, $g^{i j}$ is in fact the only
invariant symmetric (and analytic) such tensor.

Assume the statement is true for $n$, and suppose $I_{n+2}$ and $ I'_{n+2}$
satisfy the above conditions. Using the uniqueness of $I_{n}$, we have
(in shorthand - notation)
\berr
g_{12} I_{n+2}        &=& f(q-1) I_{n}  \\
g_{12} I'_{n+2} &=& f'(q-1) I_{n}
\err
where the $f(q-1)$ are analytic,
because the left - hand sides are invariant, symmetric and analytic.
Then $J_{n+2} = f' I_{n+2} - f I'_{n+2}$ is
symmetric, analytic, and satisfies $g_{12} J_{n+2} = 0$. It remains to show
that $J=0$.

Since $J$ is analytic, we can write
\beq
J^{i_1 ... i_{n}} = \sum_{k=n_0}^{\infty} (q-1)^k J_{(k)}^{i_1 ... i_n}.
\eq
$(q-1)^{-n_0} J^{i_1 ... i_n}$ has all the properties of $J$
and has a well-defined, nonzero limit as $q \rightarrow 1$; so we may assume
$J_{(0)} \neq 0$.

Now consider invariance,
\beq
J^{i_1 ... i_n} = A^{i_1}_{j_1} ... A^{i_n}_{j_n} J^{j_1 ... j_n}.
\label{J_inv}
\eq
This equation is valid for all $q$, and we can take the limit $q
\rightarrow 1$.
Then $A^i_j$ generate the commutative algebra of
functions
on the classical Lie group $O(N)$, and $J$ becomes $J_{(0)}$, which is just
a classical
tensor. Now $(P^-)^{i_l i_{l+1}}_{j_l j_{l+1}} J^{j_1 ... j_n} = 0$
implies that $J_{(0)}$ is symmetric
for neighboring indices, and therefore it is completely symmetric.
With $g_{12} J = 0$, this implies that $J_{(0)}$ is totally traceless
(classically!).
But there exists no
totally symmetric traceless invariant tensor for $O(N)$.
This proves uniqueness. In particular,  $I^{i_1 ... i_{n}} =
\l_n (\Delta^{n/2} x^{i_1} ... x^{i_{n}})$ obviously satisfies the
assumptions of the proposition; it is analytic, because in evaluating the
Laplacians, only the metric and the $\hat R$ - matrix are involved,
which are both analytic. Statement
(\ref{I_norm}) now follows because $x^2$ is central.
\endproof

Such invariant tensors have also been considered in \cite{fiore}
(where they are called S), as well as the explicit form involving the
Laplacian. Our contribution is a self - contained proof of their uniqueness.
So $<t^{i_1} ... t^{i_n}>_t \equiv I^{i_1 ... i_n}$ for even $n$ (and 0 for
odd $n$) defines the unique invariant integral over $\Sq$,
which thus coincides with the definition given in \cite{weich}.

{}From now on we only consider $N\geq 3$ since for $N=1,2$, Euclidean space
is undeformed. The following lemma is the origin of the cyclic properties
of invariant tensors. For quantum groups, the square of the antipode
is usually not 1. For $(S)O_q(N)$, it is
generated by the $D$ - matrix: $S^2 A^i_j = D^i_l A^l_k (D^{-1})^k_j$ where
$D^i_l = g^{ik} g_{lk}$ (note that $D$ also defines the quantum trace). Then

\begin{lemma}     \label{inv_lemma}
For any invariant tensor  $ J^{i_1 ... i_n} =
A^{i_1}_{j_1} ... A^{i_n}_{j_n} J^{j_1 ... j_n}$,
$D^{i_1}_{l_1} J^{i_2 ... l_1}$ is invariant too:
\beq  \label{D_cycl}
A^{i_1}_{j_1} ... A^{i_n}_{j_n} D^{j_1}_{l_1} J^{j_2 ... l_1} =
D^{i_1}_{l_1} J^{i_2 ... l_1}
\eq
\end{lemma}

\proof
{}From the above, (\ref{D_cycl}) amounts to
\beq
(S^{-2} A^{i_1}_{j_1}) A^{i_2}_{j_2} ... A^{i_n}_{j_n} J^{j_2 ... j_n j_1}
= J^{i_2 ... i_n i_1}.
\eq
Multiplying with $S^{-1}A^{i_0}_{i_1}$ from the left
and using $S^{-1}(ab) = (S^{-1}b)(S^{-1}a)$ and
$(S^{-1}A^{i_1}_{j_1}) A^{i_0}_{i_1} = \delta^{i_0}_{j_1}$, this
becomes
\beq
A^{i_2}_{j_2} ... A^{i_n}_{j_n} J^{j_2 ... j_n i_0}
= S^{-1}A^{i_0}_{i_1} J^{i_2 ... i_n i_1}.
\eq
Now multiplying with $A^{l_0}_{i_0}$ from the right, we get
\beq
A^{i_2}_{j_2} ... A^{i_n}_{j_n} A^{l_0}_{i_0} J^{j_2 ... j_n i_0}
= \delta^{l_0}_{i_1} J^{i_2 ... i_n i_1}.
\eq
But the (lhs) is just $J^{i_2 ... i_n l_0}$ by invariance and thus equal to
the (rhs).
\endproof

We can now show a number of properties of the integral over the sphere:

\begin{theorem} \label{spher_thm}
\berr
    \overline{<f(t)>_t}     \quad   &=& <\overline{f(t)}>_t
        \label{reality_t} \\
   <\overline{f(t)} f(t)>_t \quad  &\geq& 0     \label{pos_def_t} \\
   <f(t) g(t)>_t           \quad  &=&     <g(t) f(Dt)>_t     \label{cyclic_t}
\err
where $(Dt)^i = D^i_j t^j$. The last statement follows from
\beq
I^{i_1 ... i_n} = D^{i_1}_{j_1} I^{i_2 ... i_n j_1}.     \label{cyclic_I}
\eq
\end{theorem}

\proof
For (\ref{reality_t}), we have to show that
$I^{j_n ... j_1} g_{j_n i_n} ... g_{j_1 i_1} = I^{i_1 ... i_n}$.
Using the uniqueness of $I$, it is enough to show that
$I^{j_n ... j_1} g_{j_n i_n} ... g_{j_1 i_1}$ is
invariant, symmetric and normalized as $I$.
So first,
\berr
A^{i_1}_{j_1} ... A^{i_n}_{j_n}
     \(I^{k_n ... k_1} g_{k_n j_n} ... g_{k_1 j_1}\)
&=& g_{l_1 i_1} ... g_{l_n i_n} \overline{A^{l_n}_{k_n} ... A^{l_1}_{k_1}}
     I^{k_n ... k_1} \nonumber \\
&=& \overline{A^{l_n}_{k_n} ... A^{l_1}_{k_1}I^{k_n ... k_1}}
                g_{l_1 i_1} ... g_{l_n i_n}  \nonumber \\
&=& \(I^{l_n ... l_1} g_{l_n i_n} ... g_{l_1 i_1}\).
\err
We have used that $I$ is real (since $g^{ij}$ and $\hat R$ are real), and
$A^{i_1}_{j_1} g_{k_1 j_1} = g_{l_1 i_1}\overline{A^{l_1}_{k_1}}$.
The symmetry condition (\ref{I_symm}) follows from standard
compatibility conditions between $\hat R$ and $g^{ij}$,
and the fact that $\hat R$ is symmetric.
The correct normalization
can be seen easily using $g^{ij} = g_{ij}$ for $q$ - Euclidean space.

To show positive definiteness (\ref{pos_def_t}), we use the observation
made by \cite{FRT} that
\beq
t^i \rightarrow A^i_j u^j
\eq
with $u^j = u_1\delta^j_1 + u_N\delta^j_N$ is an embedding $\Sq
\rightarrow Fun(O_q(N))$ for $u_1 u_N = (q^{(N-2)/2} + q^{(2-N)/2})^{-1}$,
since $(P^-)^{ij}_{kl} u^k u^l = 0$ and $g_{ij} u^i u^j = 1$. In fact,
this embedding also respects the star - structure if one chooses
$u_N = u_1 q^{1-N/2}$ and real.
Now one can write the integral over $\Sq$ in terms of the Haar - measure
on the compact quantum group $O_q(N,\reals)$ \cite{woronowich,podles}. Namely,
\beq
<t^{i_1} ... t^{i_n}>_t
  = <A^{i_1}_{j_1} ... A^{i_n}_{j_n}>_A u^{j_1} ... u^{j_n}
\equiv <A^{\und{i}}_{\und{j}}>_A u^{\und{j}},
\eq
(in short notation)
since the Haar - measure $<>_A$ is left (and right) - invariant
$<A^{\und{i}}_{\und{j}}>_A = A^{\und{i}}_{\und{k}} <A^{\und{k}}_{\und{j}}>_A =
<A^{\und{i}}_{\und{k}}>_A A^{\und{k}}_{\und{j}}$
and analytic, and the
normalization condition is satisfied as well.
Then
$<\overline{t^{\und{i}}} t^{\und{j}}>_t =
   <\overline{A^{\und{i}}_{\und{k}}} A^{\und{j}}_{\und{r}}>_A
u^{\und{k}} u^{\und{r}}$
and for $f(t) = \sum f_{\und{i}} t^{\und{i}}$ etc.,
\berr
<\overline{f(t)} g(t)>_t &=& \overline{f_{\und{i}}} g_{\und{j}}
<\overline{A^{\und{i}}_{\und{k}}}  A^{\und{j}}_{\und{r}}>_A
  u^{\und{k}}  u^{\und{r}}
          = <\overline{\(f_{\und{i}}A^{\und{i}}_{\und{k}}u^{\und{k}}\)}
           \(g_{\und{j}}A^{\und{j}}_{\und{r}}u^{\und{r}}\)>_A \nonumber \\
 &=& <\overline{f(Au)} g(Au)>_A.
\err
This shows that the integral over $\Sq$ is positive definite, because the
Haar - measure over compact quantum groups is positive definite
\cite{woronowich}, cp. \cite{kornw}.

Finally we show the cyclic property (\ref{cyclic_I}). (\ref{cyclic_t}) then
follows immediately.
For $n=2$, the statement is obvious: $g^{ij} = D^i_k g^{jk}$.\\
Again using a shorthand - notation, define
\beq
J^{12...n} = D_1 I^{23 ... n1}.
\eq
Using the previous proposition, we only have to show that $J$ is symmetric,
invariant, analytic and properly normalized. Analyticity is obvious.
The normalization follows immediately by induction, using the property
shown in proposition (\ref{spher_int}).
Invariance of $J$ follows from the above lemma.
It remains to show that $J$ is symmetric, and the only nontrivial part of
that is $(P^-)_{12} J^{12...n} = 0$.
Define
\beq
\tilde J^{12...n} = (P^-)_{12} J^{12...n},
\eq
so $\tilde J$ is invariant, antisymmetric and
traceless in the first two indices $(12)$, symmetric in the remaining indices
(we will say that such a tensor has the ISAT property), and analytic.
It is shown below that there is no such $\tilde J$ for $q=1$ (and $N \geq 3$).
Then as in proposition (\ref{spher_int}), the leading term of the expansion of
$\tilde J$ in $(q-1)$ is classical and therefore
vanishes, which proves that $\tilde J=0$ for any $q$.

So from now on $q=1$. We show by induction that $\tilde J =0$.
This is true for $n=2$: there is no invariant antisymmetric traceless
tensor with 2 indices (for $N \geq 3$).
Now assume the statement is true for $n$ even, and
that $\tilde J^{12...(n+2)}$ has the ISAT property.
Define
\beq
K^{12...n} = g_{(n+1),(n+2)} \tilde J^{12...(n+2)}.
\eq
$K$ has the ISAT property, so by the induction assumption
\beq
K=0   \label{K}.
\eq
Define
\beq
M^{145...(n+2)} = g_{23} \tilde J^{12...(n+2)} =
              {\cal S}_{14} M^{145...(n+2)} + {\cal A}_{14} M^{145...(n+2)}
\eq
where ${\cal S}$ and ${\cal A}$ are the classical symmetrizer and
antisymmetrizer.
Again by the induction assumption, ${\cal A}_{14} M^{145...(n+2)} = 0$
(it satisfies the ISAT property). This shows that $M$ is symmetric
in the first two indices $(1,4)$. Together with the definition of $M$,
this implies
that $M$ is totally symmetric. Further, $g_{14} M^{145...(n+2)} =
g_{14} g_{23} \tilde J^{12...(n+2)} = 0$ because $\tilde J$ is
antisymmetric in $(1,2)$. But then $M$ is totally traceless, and as in
proposition (\ref{spher_int}) this implies $M=0$.
Together with (\ref{K}) and the ISAT
property of $\tilde J$, it follows that $\tilde J$ is totally traceless.
So $\tilde J$ corresponds to a certain Young tableaux, describing a
larger - than - one dimensional
irreducible representation of $O(N)$.
However, $\tilde J$ being
invariant means that it is a trivial one - dimensional representation.
This is a contradiction and proves $\tilde J=0$.

\endproof

Property (\ref{pos_def_t})\footnote{as was pointed out to me by G. Fiore,
positivity is also implied by results in \cite{fiore}}
in particular means that one can now define the
Hilbertspace of square - integrable functions on $\Sq$. The same will be true
for the integral on the entire Quantum Euclidean space.

The cyclic property (\ref{cyclic_t}) is a strong constraint on
$I^{i_1 ... i_n}$ and could actually be used to
calculate it recursively, besides its obvious interest in its own.

\section{Integral over quantum Euclidean space}

It is now easy to define an integral over quantum Euclidean space.
Since the invariant length $r^2= g_{ij} x^i  x^j$ is central, we can use its
square root $r$ as well, and write any function
on quantum Euclidean space in the form $f(x^i) = f(t^i,r)$.
We then define its integral
to be
\beq
<f(x)>_x = <<f(t,r)>_t(r) \cdot r^{N-1}>_r,    \label{integral_x}
\eq
where $<f(t,r)>_t(r)$ is a classical, analytic function in $r$,
and $<g(r)>_r$ is some
linear functional in $r$, to be determined by requiring Stokes theorem.
It is essential that this radial integral $<g(r)>_r$ is really a functional
of the
{\em analytic continuation of g(r)} to a function on the (positive) real line.
Only then one obtains a large class of integrable functions,
and this concept of integration over the entire space agrees
with the classical one. This is
also the reason why the Gaussian integration procedure suggested e.g.
in \cite{fiore,KM}
works only for a very small class of functions of the form $p(x) g_{\a}(x)$
where $g_{\a}(x)$ is a Gaussian
and $p(x)$ is a certain class of power - series. It diverges as soon as $p(x)$
has a finite radius of "convergence",
such as e.g. $1/(r^2+1)$; even classically,
one cannot do such integrals term by term. This problem
does not occur for the spherical integral.

It will turn out that Stokes theorem e.g. in the form $<\dl_i f(x)>_x = 0$
holds if and only if
the radial integral satisfies the scaling property
\beq
<g(qr)>_r = q^{-1} <g(r)>_r.   \label{scaling_r}
\eq
This can be shown directly; we will instead give a more elegant proof later.
This scaling property is obviously satisfied by an arbitrary
superposition of Jackson - sums,
\beq
<f(r)>_r = \int_1^q dr_0 \mu(r_0)\sum_{n=-\infty}^{\infty} f(q^nr_0) q^{n}
                                        \label{jackson_int}
\eq
with arbitrary (positive) "weight" $\mu(r_0) >0$. If $\mu(r_0)$ is a
delta - function,
this is simply a Jackson - sum; for $\mu(r_0) = 1$, one obtains
the classical radial integration
\beq
<f(r) r^{N-1}>_r = \int_1^q dr_0 q^n  (q^n r_0)^{N-1}
\sum_{n=-\infty}^{\infty} f(q^nr_0) = \int_0^{\infty} dr r^{N-1}f(r).
\eq
For Gaussian - integrable functions, all of these uncountable
choices are equivalent (and of course agree with that definition).
This is in general not true for
functions integrable in this radial sense,
i.e. for which the above is finite, and
shows again that this class is indeed larger.
The classical integral over $r$ however is the unique choice for which the
scaling property (\ref{scaling_r}) holds for any positive real number,
not just for powers of $q$.

The properties of the integral over $\Sq$  generalize immediately
to the Euclidean case:
\begin{theorem}
\berr
 \overline{<f(x)>_x}    \quad  &=& <\overline{f(x)}>_x  \label{reality_x} \\
 <\overline{f(x)} f(x)>_x \quad &\geq& 0    \label{pos_def_x} \\
 <f(x) g(x)>_x             \quad  &=&  <g(x) f(Dx)>_x,     \label{cyclic_x}
\err
and
\beq
<f(qx)>_x \quad   = q^{-N} <f(x)>_x     \label{scaling_x}
\eq
if and only if (\ref{scaling_r}) holds.
\end{theorem}

\proof
Immediately from theorem (\ref{spher_thm}), (\ref{scaling_r}) and
(\ref{integral_x}), using $Dr = r$ and $\mu(r_0)>0$.
\endproof

(\ref{reality_x}) and (\ref{scaling_x}) were already known for the
special case of the Gaussian integral \cite{fiore}\footnote{
it was pointed out to me by G. Fiore that positivity for certain
classes of functions was shown in \cite{fiore_thesis}}.

\section{Integration of forms}

It turns out to be very useful to consider not only integrals over functions,
but also over forms, just like classically. As was mentionned before,
there exists a unique
$N$ - form $dx^{i_1} ... dx^{i_N} = \eps_q^{i_1 ... i_N} \d x$, and we define
\beq    \label{int_x_form}
\int_x \d x f(x) = <f(x)>_x,
\eq
i.e. we first commute $\d x$ to the left, and then take
the integral over the function on the right.
Then the two statements of Stokes theorem $<\dl_i f(x)>_x = 0$ and
$\int_x d\omega_{N-1} = 0$ are equivalent.

The following observation by Bruno Zumino \cite{zumino}
will be very useful:
there is a one - form
\beq
\om = \frac {q^2}{(q+1)r^2} d(r^2) = q\frac 1r dr = dr \frac 1{r}
\eq
where $r dx^i = q dx^i r$,
which generates the calculus on quantum Euclidean space by
\beq
[\om, f]_{\pm} = (1-q) df
\eq
for any form $f$ with the appropriate grading. It satisfies
\beq
d\om = \om^2 = 0.
\eq
We define the integral of a $N$ - form over the sphere $r\cdot\Sq$
with radius $r$ by

\beq
\int\limits_{r\cdot\Sq} \d x f(x) = \om r^N <f(x)>_t = dr r^{N-1} <f(x)>_t,
            \label{int_S_Nform}
\eq
which is a one - form in $r$, as classically. It satisfies
\beq
\int\limits_{r\cdot\Sq} q^N\d x f(qx) = \int\limits_{qr\cdot\Sq} \d x f(x)
         \label{scaling_S_Nform}
\eq
where $(dr f(r))(qr) = q dr f(qr)$.
Now defining $\int_r dr g(r) = <g(r)>_r$,
(\ref{int_x_form}) can be written as
\beq
\int_x \d x f(x) = \int_r (\int\limits_{r\cdot\Sq} \d x f(x)).
\eq
The scaling property (\ref{scaling_r}), i.e.
$\int_x \d x f(qx) = q^{-N}\int_x \d x f(x)$
holds if and only if the radial integrals satisfies
\beq   \label{scaling_r_forms}
\int_r dr f(qr) = q^{-1} \int_r dr f(r).
\eq
We can also define the integral
of a $(N-1)$ form $\a_{N-1}(x)$ over the sphere with radius $r$:
\beq   \label{int_sphere_N-1}
\int\limits_{r\cdot\Sq} \a_{N-1}
= \om^{-1} (\int\limits_{r\cdot\Sq} \om \a_{N-1}).
\eq
The $\om^{-1}$ simply cancels the explicit $\om$ in (\ref{int_S_Nform}),
and it reduces to the correct classical limit for $q=1$.

The epsilon - tensor satisfies  the cylic property:

\begin{prop}
\beq   \label{cycl_eps}
\eps_q^{i_1 ... i_N} = (-1)^{N-1} D^{i_i}_{j_1} \eps_q^{i_2 ... i_N j_1}.
\eq
\end{prop}

\proof
Define
\beq
\kappa^{12...N} = (-1)^{N-1} D^1 \eps_q^{23...N1}
\eq
in shorthand - notation again. Lemma (\ref{inv_lemma}) shows that
$\kappa$ is invariant. $\kappa^{12...N}$ is traceless and
($q$ -) antisymmetric in (23...N).
Now $g_{12} \kappa^{12...N}=0$ because
there exists no invariant, totally antisymmetric traceless tensor
with $(N-2)$ indices for $q=1$,
so by analyticity there is none for arbitrary $q$. Similarly from
the theory of irreducible representations of $SO(N)$ \cite{weyl},
${P^+}_{12} \kappa^{12...N} = 0$
where ${P^+}$ is the $q$ - symmetrizer, $1=P^+ + P^- + P^0$.
Therefore $\kappa^{12...N}$ is totally
antisymmetric and traceless (for neighboring indices), invariant
and analytic.
But there exists only one such tensor up to normalization
(which can be proved similarly), so
$\kappa^{12...N} = f(q) \eps_q^{12...N}$. It remains to show $f(q) =1$.
By repeating the above, one gets
$\eps_q^{12...N} = (f(q))^N (\det D) \eps_q^{12...N}$
(here $12...N$ stands for the {\em numbers} 1,2,...,$N$),
and since $\det D = 1$, it follows $f(q) = 1$ (times a $N$-th root of unity,
which is fixed by the classical limit).
\endproof

Now consider a $k$ - form
$\a_k(x) = dx^{i_1} .... dx^{i_k} f_{i_1 ... i_k} (x)$
and a $(N-k)$ - form $\b_{N-k}(x)$.
Then the following cyclic property for the integral over forms holds:
\begin{theorem} \label{cyclic_forms}

\beq   \label{cyclic_forms_any}
\int\limits_{r\cdot\Sq} \a_k(x) \b_{N-k}(x) =
        (-1)^{k(N-k)} \int\limits_{q^{-k}r\cdot\Sq} \b_{N-k}(x) \a_k(q^N Dx)
\eq
where
$\a_k(q^N Dx) = (q^N Ddx)^{i_1} .... (q^N Ddx)^{i_k} f_{i_1 ... i_k} (q^N Dx)$.

In particular, when $\a_k$ and $\b_{N-k}$ are forms on $\Sq$,
i.e. they involve only $dx^i \frac 1r$ and $t^i$,
then
\beq    \label{cyclic_forms_spher}
\int\limits_{\Sq} \a_k(t) \b_{N-k}(t)
= (-1)^{k(N-k)}\int\limits_{\Sq} \b_{N-k}(t) \a_k(Dt).
\eq
On Euclidean space,
\beq    \label{cyclic_forms_x}
\int_x \a_k(x) \b_{N-k}(x) = (-1)^{k(N-k)}\int_x \b_{N-k}(x) \a_k(q^N Dx)
\eq
if and only if (\ref{scaling_r_forms}) holds.
\end{theorem}
Notice that on the sphere, $\d x f(t) = f(t)\d x$.

\proof
We only have to show that
\beq   \label{cycl_spher_f}
\int\limits_{r\cdot\Sq} f(x) \d x g(x) =
   \int\limits_{r\cdot\Sq} \d x g(x) f(q^N Dx)
\eq
and
\beq   \label{cycl_spher_diff}
\int\limits_{r\cdot\Sq} dx^i \b_{N-1}(x) =
  (-1)^{N-1} \int\limits_{q^{-1}r\cdot\Sq} \b_{N-1}(x) (q^N Ddx)^i.
\eq
(\ref{cycl_spher_f}) follows immediately from (\ref{cyclic_t}) and
$x^i \d x = \d x q^N x^i$.

To see (\ref{cycl_spher_diff}), we can assume that $\b_{N-1}(x)
= dx^{i_2} ... dx^{i_N} f(x)$.
The commutation relations
$x^i dx^j = q\hat R^{ij}_{kl} dx^k x^l$
are equivalent to
\berr   \label{CR_compact}
f(q^{-1}x) dx^j &=& \R((dx^j)_{(a)} \tens f_{(1)}) (dx^j)_{(b)}
           (f(x))_{(2)} \nonumber\\
          &=& (dx^j \triangleleft \R^1) (f(x)\triangleleft \R^2)
\err
where $\R = \R^1 \tens \R^2$ is the universal $\R$ for $SO_q(N)$,
using its quasitriangular property
and $\R(A^j_k \tens A^i_l) = \hat R^{ij}_{kl}$.
$f\triangleleft Y = <Y,f_{(1)}> f_{(2)}$ is the right action induced by
the left coaction (\ref{left_coaction}) of an element
$Y \in {\cal U}_q(SO(N))$.
Now invariance of the integral implies
\beq  \label{int_inv}
(dx^j \triangleleft \R^1) <f(x)\triangleleft \R^2>_t = dx^j <f(x)>_t,
\eq
because $\R^1 \tens\eps(\R^2) = 1$.
Using this, (\ref{CR_compact}), (\ref{scaling_S_Nform}) and
(\ref{int_S_Nform}), the (rhs) of (\ref{cycl_spher_diff}) becomes
\berr
(-1)^{N-1} \int\limits_{q^{-1}r\cdot\Sq} \b_{N-1}(x) q^N D^i_j dx^j
     &=& (-1)^{N-1} D^i_j \int\limits_{r\cdot\Sq} dx^{i_2} ... dx^{i_N}
           f(q^{-1}x) dx^j \nonumber\\
     &=& (-1)^{N-1} D^i_j \eps^{i_2 ... i_N j} \om r^N <f(x)>_t  \nonumber\\
     &=& \eps^{i i_2 ... i_N} \om r^N <f(x)>_t \nonumber\\
     &=& \int\limits_{r\cdot\Sq} dx^i \b_{N-1}(x),
\err
using (\ref{cycl_eps}). This shows (\ref{cycl_spher_diff}), and
(\ref{cyclic_forms_spher}) follows immediately. (\ref{cyclic_forms_x})
then follows from (\ref{scaling_r_forms}).

\endproof

Another way to show (\ref{cycl_spher_diff}) following an idea of
Branislav Jurco \cite{jurco} is to use
\beq
\int\limits_{r\cdot\Sq} (\a_k \triangleleft SY)\b_{N-k}
= \int\limits_{r\cdot\Sq}  \a_k (\b_{N-k} \triangleleft Y)
\eq
to move the action of $\R^2$ in (\ref{CR_compact}) to the left picking up
$\R^1 S \R^2$, which generates
the inverse square of the antipode and thus corresponds to
the $D^{-1}$ - matrix. This approach however cannot
show (\ref{cyclic_t}) or (\ref{cyclic_x}),
because the commutation relations of functions are more complicated.

(\ref{cyclic_forms_any}) shows in particular that the definition
(\ref{int_sphere_N-1}) is natural,
i.e. it essentially does not matter on which side one multiplies with $\om$.
Now we immediately obtain Stokes theorem for the integral over
quantum Euclidean space, if and only if (\ref{scaling_r_forms}) holds.
Noticing that $\om(q^NDx) = \om(x)$, (\ref{cyclic_forms_x}) implies
\berr
\int_x d\a_{N-1}(x) &=& \frac{1}{1-q} \int_x [\om, \a_{N-1}]_{\pm}
               \nonumber \\
       &\propto& \int_x \om \a_{N-1} - (-1)^{N-1} \a_{N-1}\om   \nonumber \\
       &=& \int_x (-1)^{N-1} \a_{N-1}\om - (-1)^{N-1} \a_{N-1}\om = 0
\err
On the sphere, we get as easily
\berr
\int\limits_{\Sq} d\a_{N-2}(t) &\propto& \int\limits_{\Sq}
                   [\om, \a_{N-2}]_{\pm}   \nonumber \\
      &=& \om^{-1} \int\limits_{\Sq} \om
               (\om \a_{N-2} - (-1)^{N-2} \a_{N-2}\om))
     =0
\err
using (\ref{cyclic_forms_spher}) and $\om^2 = 0$.

It is remarkable that these simple proofs only work for $q\neq 1$,
nevertheless the statements reduce to the classical Stokes theorem
for $q \rightarrow 1$.
This shows the power of the $q$ - deformation technique.

One can actually obtain a version of Stokes theorem
with spherical boundary terms.
Define
\beq
\int\limits_{q^k r_0}^{q^l r_0} \om f(r) =
\int\limits_{q^k r_0}^{q^l r_0} dr \frac 1r f(r) =
 (q-1) \sum_{n=k}^{l-1} f(r_0 q^n),  \label{int_r_finite}
\eq
which reduces to the correct classical limit, because the
(rhs) is a Riemann sum.
Define
\beq
\int\limits_{q^k r_0\cdot\Sq}^{q^l r_0\cdot\Sq} \a_N(x)
= \int\limits_{q^k r_0}^{q^l r_0}
(\int\limits_{r \cdot\Sq}\a_N(x)),
\eq
For $l \rightarrow \infty$ and  $k \rightarrow -\infty$, this becomes
an integral over Euclidean space as defined before.
Then
\berr
\int\limits_{q^k r_0\cdot\Sq}^{q^l r_0\cdot\Sq} d\a_{N-1}
    &=& \frac{1}{1-q}\int\limits_{q^k r_0}^{q^l r_0}
    (\int\limits_{r\cdot\Sq}\om \a_{N-1} -(-1)^{N-1} \a_{N-1}\om) \nonumber\\
    &=& \frac{1}{1-q} \int\limits_{q^k r_0}^{q^l r_0}
       \big(\int\limits_{r \cdot\Sq}\om \a_{N-1}
           - \int\limits_{qr \cdot\Sq}\om\a_{N-1} )     \nonumber\\
    &=& \int\limits_{q^l r_0 \cdot\Sq}\a_{N-1} -
        \int\limits_{q^k r_0 \cdot\Sq}\a_{N-1}.
\err
In the last line, (\ref{int_S_Nform}), (\ref{int_sphere_N-1}) and
(\ref{int_r_finite}) was used.

\section{Acknowledgements}
It is a pleasure to thank Prof. Bruno Zumino for many useful
discussions and encouragement. I also wish to thank Chris Chryssomalakos,
Chong -Sun Chu, Pei - Ming Ho, Branislav Jurco and Bogdan Morariu for
very useful discussions, and Gaetano Fiore for pointing out his results
to me.
This work was supported in part by the Director, Office of
Energy Research, Office of High Energy and Nuclear Physics, Division of
High Energy Physics of the U.S. Department of Energy under Contract
DE-AC03-76SF00098, in part by the National Science Foundation under
grant PHY-90-21139 and in part by a Regents fellowship of the
University of California.

\end{document}